\documentclass[12pt]{article}

\usepackage{amsfonts}

\usepackage{amsmath,amssymb}

 \numberwithin{equation}{subsection}

\begin{document}

\title{Uniqueness of the representation for $G$-martingales with finite variation}
\author{Yongsheng Song
\footnote{ Y. Song was  supported by the National Basic Research
Program of China (973 Program) (No.2007CB814902), Key Lab of Random
Complex Structures and Data Science, Chinese Academy of Sciences
(Grant No. 2008DP173182).
}\\
\small Academy of Mathematics and Systems Science, \\
\small Chinese Academy of Sciences, Beijing, China;\\
\small yssong@amss.ac.cn}

\date{}

\maketitle
\begin{abstract}

Our purpose is to prove the uniqueness  of the representation for
$G$-martingales with finite variation.

\end{abstract}


{\bf MSC-classification}:  60G48, 60G44
\maketitle
\section{Introduction }

 \ \ \  \ In [P07b], processes in form of $\int_0^t\eta_sd\langle
 B\rangle_s-\int_0^t2G(\eta_s)ds$, $\eta\in M^1_G(0,T)$ are proved
 to be
 $G$-martingales. However, the uniqueness of the representation
 remains unresolved. In order to prove the uniqueness, we must
 find ways to distinguish the two classes of processes in forms of $\int_0^t\eta_sd\langle
 B\rangle_s$ and $\int_0^t\zeta_sds$, $\eta, \zeta\in M^1_G(0,T)$.

 For a process $\{K_t\}$ with finite variation, motivated by [Song10], we define
$$d(K):=\limsup_{n\rightarrow\infty}\hat{E}[\int_0^T\delta_n(s)dK_s],$$
 where, for $n\in N$, $\delta_n(s)$ is defined in the following way:
$$\delta_n(s)=\sum_{i=0}^{n-1}(-1)^i1_{]\frac{iT}{n}, \frac{(i+1)T}{n}]}(s), \textmd{\ for \ all} \ s\in[0,T].$$
We prove that $d(K)=0$ if $K_t=\int_0^t\zeta_sds$ for some $\zeta\in
M^1_G(0,T)$ and that $d(K)>0$ if $K_t=\int_0^t\eta_sd\langle
 B\rangle_s$ for some
$\eta\in M^1_G(0,T)$ such that $\hat{E}[\int_0^T|\eta_s|ds]>0$. By
this, we distinguish these two classes of processes completely:

If $\int_0^t\eta_sd\langle
 B\rangle_s=\int_0^t\zeta_sds$, for some $\eta, \zeta\in
 M^1_G(0,T)$, then we have
 $$\hat{E}[\int_0^T|\eta_s|ds]=\hat{E}[\int_0^T|\zeta_s|ds]=0.$$

As an application, we prove the uniqueness  of the representation
for $G$-martingales with finite variation.

This article is organized as follows: In section 2, we recall some
basic notions and results of $G$-expectation and the related space
of random variables. In section 3, we present the main results and
some corollaries. In section 4, we give the proofs to the main
results.

\section{Preliminaries }
We recall some basic notions and results of $G$-expectation and the
related space of random variables. More details of this section can
be found in [P07a, P07b, P08, P10].

\noindent {\bf Definition 2.1} Let $\Omega$
 be a given set and let ${\cal H}$ be a linear space of real valued
functions defined on $\Omega$
 with $c \in {\cal H}$ for all constants $c$. ${\cal H}$ is considered as the
space of  ¡°random variables¡±. A sublinear expectation $\hat{E}$ on
${\cal H}$ is a functional $\hat{E}: {\cal H}\rightarrow R $
satisfying the following properties: for all $X, Y \in {\cal H}$, we
have

(a) Monotonicity: If $X\geq Y$ then $\hat{E}(X) \geq \hat{E} (Y)$.

(b) Constant preserving: $\hat{E} (c) = c$.

(c) Sub-additivity: $\hat{E}(X)-\hat{E}(Y) \leq \hat{E}(X-Y)$.

(d) Positive homogeneity: $\hat{E} (\lambda X) = \lambda \hat{E}
(X)$, $\lambda \geq 0$.

\noindent$(\Omega, {\cal H}, \hat{E})$ is called a sublinear
expectation space.

\noindent {\bf Definition 2.2} Let $X_1$ and $X_2$ be two
$n$-dimensional random vectors defined respectively in sublinear
expectation spaces $(\Omega_1, {\cal H}_1, \hat{E}_1)$ and
$(\Omega_2, {\cal H}_2, \hat{E}_2)$. They are called identically
distributed, denoted by $X_1 \sim X_2$, if $\hat{E}_1[\varphi(X_1)]
= \hat{E}_2[\varphi(X_2)]$, $\forall \varphi\in C_{l, Lip}(R^n)$,
where $ C_{l, Lip}(R^n)$ is the space of real continuous functions
defined on $R^n$ such that $$|\varphi(x) - \varphi(y)| \leq C(1 +
|x|^k + |y|^k)|x - y|, \forall x, y \in R^n,$$ where $k$ depends
only on $\varphi$.

\noindent {\bf Definition 2.3} In a sublinear expectation space
$(\Omega, {\cal H}, \hat{E})$ a random vector $Y = (Y_1,
\cdot\cdot\cdot, Y_n)$, $Y_i \in {\cal H}$ is said to be independent
to another random vector $X = (X_1, \cdot\cdot\cdot, X_m)$, $X_i \in
{\cal H}$ under $\hat{E}(\cdot)$, denoted by $Y\bot X$, if for each
test function $\varphi\in C_{l, Lip}(R^m\times R^n)$ we have
$\hat{E}[\varphi(X, Y )] = \hat{E}[\hat{E} [\varphi(x, Y )]_{x=X}]$.

\noindent {\bf Definition 2.4} ($G$-normal distribution) A
d-dimensional random vector $X = (X_1, \cdot\cdot\cdot,X_d)$ in a
sublinear expectation space $(\Omega, {\cal H}, \hat{E})$ is called
$G$-normal distributed if for each $a, b\in R$ we have $$aX +
b\hat{X}\sim \sqrt{a^2 + b^2}X,$$  where $\hat{X}$ is an independent
copy of $X$. Here the letter $G$ denotes the function $$G(A) :=
\frac{1 }{2}\hat{ E}[(AX,X)] : S_d \rightarrow R,$$  where $S_d$
denotes the collection of $d\times d$ symmetric matrices.

The function $G(\cdot) : S_d \rightarrow R$ is a monotonic,
sublinear mapping on $S_d$ and $G(A) = \frac{1 }{2}\hat{
E}[(AX,X)]\leq \frac{1 }{2}|A|\hat{ E}[|X|^2]=:\frac{1
}{2}|A|\bar{\sigma}^2$ implies that there exists a bounded, convex
and closed subset $\Gamma\subset S_d^+$ such that
\begin {eqnarray}\label{eqn1}
G(A)=\frac{1 }{2}\sup_{\gamma\in \Gamma}Tr(\gamma A).
\end {eqnarray}
If there exists some $\beta>0$ such that $G(A)-G(B)\geq \beta
Tr(A-B)$ for any $A\geq B$, we call the $G$-normal distribution is
non-degenerate.

\noindent {\bf Definition 2.5} i) Let $\Omega_T=C_0([0, T]; R^d)$
with the supremum norm, $ {\cal H}^0_T:=\{\varphi(B_{t_1},...,
B_{t_n})| \forall n\geq1, t_1, ..., t_n \in [0, T], \forall \varphi
\in C_{l, Lip}(R^{d\times n})\}$, $G$-expectation is a sublinear
expectation defined by
$$\hat{E}[\varphi( B_{t_1}-B_{t_0}, B_{t_2}-B_{t_1} ,
\cdot\cdot\cdot, B_{t_m}- B_{t_{m-1}} )]$$$$ = \tilde{E}
[\varphi(\sqrt{t_1-t_0}\xi_1, \cdot\cdot\cdot, \sqrt{t_m
-t_{m-1}}\xi_m)],$$ for all $X=\varphi( B_{t_1}-B_{t_0},
B_{t_2}-B_{t_1} , \cdot\cdot\cdot, B_{t_m}- B_{t_{m-1}} )$, where
$\xi_1, \cdot\cdot\cdot, \xi_n$ are identically distributed
$d$-dimensional $G$-normal distributed random vectors in a sublinear
expectation space $(\tilde{\Omega}, \tilde{\cal H},\tilde{ E})$ such
that  $\xi_{i+1}$ is independent to $(\xi_1, \cdot\cdot\cdot,
\xi_i)$ for each $i = 1, \cdot\cdot\cdot,m$. $(\Omega_T, {\cal
H}^0_T, \hat{E})$ is called a $G$-expectation space.

ii) For $t\in [0, T]$ and $\xi=\varphi(B_{t_1},..., B_{t_n})\in
{\cal H}^0_T$, the conditional expectation defined by(there is no
loss of generality, we assume $t=t_i$) $$\hat{E}_{t_i}[\varphi(
B_{t_1}-B_{t_0}, B_{t_2}-B_{t_1} , \cdot\cdot\cdot, B_{t_m}-
B_{t_{m-1}} )]$$$$=\tilde{\varphi}( B_{t_1}-B_{t_0}, B_{t_2}-B_{t_1}
, \cdot\cdot\cdot, B_{t_i}- B_{t_{i-1}} ),$$ where
$$\tilde{\varphi}(x_1, \cdot\cdot\cdot, x_i)=\hat{E}[\varphi( x_1,
\cdot\cdot\cdot,x_i, B_{t_{i+1}}- B_{t_{i}}, \cdot\cdot\cdot,
B_{t_m}- B_{t_{m-1}} )].$$

Define $\|\xi\|_{p, G}=[\hat{E}(|\xi|^p)]^{1/p}$ for $\xi\in{\cal
H}^0_T$ and $p\geq1$. Then $\forall t\in[0, T]$, $\hat{E}_t(\cdot)$
is a continuous mapping on ${\cal H}^0_T$ with norm $\|\cdot\|_{1,
G}$ and therefore can be extended continuously to the completion
$L^1_G(\Omega_T)$ of ${\cal H}^0_T$ under norm $\|\cdot\|_{1, G}$.

Let $L_{ip}(\Omega_T):=\{\varphi(B_{t_1},..., B_{t_n})|  n\geq1,
t_1, ..., t_n \in [0, T],  \varphi \in C_{b, Lip}(R^{d\times n})\},$
where $C_{b, Lip}(R^{d\times n})$ denotes the set of bounded
Lipschitz functions on $R^{d\times n}$. [DHP08] proved that the
completions of $C_b(\Omega_T)$, ${\cal H}^0_T$ and
$L_{ip}(\Omega_T)$ under $\|\cdot\|_{p,G}$ are the same and we
denote them by $L^p_G(\Omega_T)$.

\noindent {\bf Definition 2.5} Let $M^0_G(0, T)$ be the collection
of processes in the following form: for a given partition $\{t_0,
\cdot\cdot\cdot, t_N\} = \pi_T$ of $[0, T]$,  $$ \eta_t(\omega) =
\sum^{N-1}_{j=0} \xi_j(\omega)1_{]t_j ,t_{j+1}]}(t),$$ where
$\xi_i\in L_{ip}(\Omega_ {t_i})$, $i = 0, 1, 2, \cdot\cdot\cdot,
N-1$. For $p\geq1$ and $\eta\in M^0_G(0, T)$, let
$\|\eta\|_{M^{p}_G}=\{\hat{E}(\int_0^T|\eta_s|^pds)\}^{1/p}$ and
denote by $M^p_G(0,T)$ the completion of $M^0_G(0, T)$ under the
norm $\|\cdot\|_{M^{p}_G}$.

\noindent {\bf Theorem 2.6}([DHP08]) There exists a tight subset
${\cal P}\subset {\cal M}_1(\Omega_T)$ such that
$$\hat{E}(\xi)=\max_{P\in {\cal P}}E_P(\xi) \ \ \textrm{for \
all} \ \xi\in{\cal H}^0_T.$$ ${\cal P}$ is called a set that
represents $\hat{E}$.

\noindent {\bf Remark 2.7} Let $(\Omega^0, \{{\cal F}^0_t\}, {\cal
F}^0, P^0 )$ be a filtered probability space and $\{W_t\}$ be a
d-dimensional Brownian motion under $P^0$.  [DHP08] proved that
$${\cal P}_M:=\{P_0\circ X^{-1}| X_t=\int_0^th_sdW_s, h\in
L^2_{\cal F}([0,T]; \Gamma^{1/2}) \}$$ is a set that represents
$\hat{E}$,  where $\Gamma^{1/2}:=\{\gamma^{1/2}| \gamma\in \Gamma\}$
and $\Gamma$ is the set in the representation of $G(\cdot)$ in the
formula (\ref {eqn1}).

\section{Main results}

In the sequel, we only consider the $G$-expectation space
$(\Omega_T, L^1_G(\Omega_T, \hat{E}))$ with $\Omega_T=C_0([0,T], R)$
and
$\overline{\sigma}^2=\hat{E}(B_1^2)>-\hat{E}(-B_1^2)=\underline{\sigma}^2\geq0$.

\noindent {\bf Proposition 3.1} For each $\eta\in M^1_G(0,T)$, let
$$d(\eta)=\limsup_{n\rightarrow\infty}\hat{E}[\int_0^T\delta_n(s)\eta_sd\langle
B\rangle_s].$$ Then
\begin {eqnarray}\label{eqn2}-\frac{\overline{\sigma}^2-\underline{\sigma}^2}{2}\hat{E}[-\int_0^T|\eta_s|ds]\leq
d(\eta)
\leq\frac{\overline{\sigma}^2-\underline{\sigma}^2}{2}\hat{E}[\int_0^T|\eta_s|ds].\end
{eqnarray}

{\bf Proof.} It suffices to prove the conclusion for $\eta\in
M^0_G(0,T)$. Let $\eta_s=\sum_{i=0}^{m-1}\xi_{t_i}1_{]t_i,
t_{i+1}]}(s)$, $\xi_{t_i}\in L^1_G(\Omega_{t_i})$,
$i=0,\cdot\cdot\cdot,m-1$.
\begin {eqnarray*}
& &\hat{E}[\int_0^T\delta_n(s)\eta_sd\langle
B\rangle_s]-\frac{\overline{\sigma}^2-\underline{\sigma}^2}{2}\hat{E}[\int_0^T|\eta_s|ds]\\
&=&\hat{E}[\sum_{i=0}^{m-1}|\xi_{t_i}|\int_{t_i}^{t_{i+1}}\delta_n(s){\textmd{sgn}(\xi_{t_i})}d\langle
B\rangle_s]-\hat{E}[\sum_{i=0}^{m-1}|\xi_{t_i}|\int_{t_i}^{t_{i+1}}\frac{\overline{\sigma}^2-\underline{\sigma}^2}{2}
ds]\\
&\leq&\sum_{i=0}^{m-1}\hat{E}[|\xi_{t_i}|(\int_{t_i}^{t_{i+1}}\delta_n(s){\textmd{sgn}(\xi_{t_i})}d\langle
B\rangle_s-\int_{t_i}^{t_{i+1}}\frac{\overline{\sigma}^2-\underline{\sigma}^2}{2}
ds)]\rightarrow0
\end {eqnarray*} as $n$ goes to infinity. So $$d(\eta)
\leq\frac{\overline{\sigma}^2-\underline{\sigma}^2}{2}\hat{E}[\int_0^T|\eta_s|ds].$$

On the other hand,
\begin {eqnarray*} & &\hat{E}[\int_0^T\delta_n(s)\eta_sd\langle
B\rangle_s]+\frac{\overline{\sigma}^2-\underline{\sigma}^2}{2}\hat{E}[-\int_0^T|\eta_s|ds]\\
&=&\hat{E}[\sum_{i=0}^{m-1}|\xi_{t_i}|\int_{t_i}^{t_{i+1}}\delta_n(s){\textmd{sgn}(\xi_{t_i})}d\langle
B\rangle_s]+\hat{E}[\sum_{i=0}^{m-1}(-|\xi_{t_i}|)\int_{t_i}^{t_{i+1}}\frac{\overline{\sigma}^2-\underline{\sigma}^2}{2}
ds]\\
&\geq&\hat{E}[\sum_{i=0}^{m-1}|\xi_{t_i}|(\int_{t_i}^{t_{i+1}}\delta_n(s){\textmd{sgn}(\xi_{t_i})}d\langle
B\rangle_s-\int_{t_i}^{t_{i+1}}\frac{\overline{\sigma}^2-\underline{\sigma}^2}{2}
ds)]\\
&\geq&\sum_{i=0}^{m-1}[-\hat{E}(|\xi_{t_i}|)a_i(n)],
\end {eqnarray*} where $a_i(n)=\max\{|\hat{E}(\int_{t_i}^{t_{i+1}}\delta_n(s)d\langle
B\rangle_s-\int_{t_i}^{t_{i+1}}\frac{\overline{\sigma}^2-\underline{\sigma}^2}{2}ds)|,
|\hat{E}(-\int_{t_i}^{t_{i+1}}\delta_n(s)d\langle
B\rangle_s-\int_{t_i}^{t_{i+1}}\frac{\overline{\sigma}^2-\underline{\sigma}^2}{2}ds)|\}\rightarrow0$
as $n$ goes to infinity. So
$$-\frac{\overline{\sigma}^2-\underline{\sigma}^2}{2}\hat{E}[-\int_0^T|\eta_s|ds]\leq
d(\eta).$$ $\Box$

\noindent {\bf Remark 3.2} (i) A straightforward corollary of
Proposition 3.1 is that if  $\int_0^T|\eta_s|ds$ is symmetric (i.e.,
$\hat{E}[\int_0^T|\eta_s|ds]=-\hat{E}[-\int_0^T|\eta_s|ds]$), the
equality $d(\eta)
=\frac{\overline{\sigma}^2-\underline{\sigma}^2}{2}\hat{E}[\int_0^T|\eta_s|ds]$
holds.

(ii) By Lemma 3.1, we could not conclude that $d(\eta)>0$ whenever
$\hat{E}[\int_0^T|\eta_s|ds]>0$, which is the conclusion of Theorem
3.3 below.

(iii) The inequalities in (\ref {eqn2}) may be strict:

Let $\eta_s=\langle B\rangle_{T/2}1_{]T/2, T]}(s)+a1_{[0,T/2]}(s)$,
$a=T(\overline{\sigma}^2-\underline{\sigma}^2)/4$.

Then
$$d(\eta)=\lim_{n\rightarrow\infty}\hat{E}[\int_0^T\delta_{2n}(s)\eta_sd\langle
B\rangle_s]=a\overline{\sigma}^2T/2,$$
$$\frac{\overline{\sigma}^2-\underline{\sigma}^2}{2}\hat{E}[\int_0^T|\eta_s|ds]=a^2+a\overline{\sigma}^2T/2,$$
$$-\frac{\overline{\sigma}^2-\underline{\sigma}^2}{2}\hat{E}[-\int_0^T|\eta_s|ds]=-a^2+a\overline{\sigma}^2T/2.$$
$\Box$

Now, we shall state the main result of this article, whose proof is
postponed to Section 4.

\noindent {\bf Theorem 3.3} For $\eta\in M^1_G(0,T)$ with
$\hat{E}[\int_0^T|\eta_s|ds]>0,$ we have
$$d(\eta)=\limsup_{n\rightarrow\infty}\hat{E}[\int_0^T\delta_n(s)\eta_sd\langle
B\rangle_s]>0.$$

\noindent {\bf Theorem 3.4} Let $\eta\in M^1_G(0,T)$. Then
$\lim_{n\rightarrow\infty}\hat{E}[\int_0^T\delta_n(s)\eta_sds]=0$.

{\bf Proof.} For $\eta\in M^0_G(0,T)$, the claim is obvious. For
$\eta\in M^1_G(0,T)$, there exists a sequence of $\{\eta^m\}\subset
M^0_G(0,T)$ such that
$\hat{E}[\int_0^T|\eta^m_s-\eta_s|ds]\rightarrow0$. Then
$|\hat{E}[\int_0^T\delta_n(s)\eta_sds]|\leq|\hat{E}[\int_0^T\delta_n(s)\eta^m_sds]|+\hat{E}[\int_0^T|\eta^m_s-\eta_s|ds]$.
First let $n\rightarrow\infty$, then let $m\rightarrow\infty$, and
we get the desired result. $\Box$

\noindent {\bf Remark 3.5} Let $(\Omega, F, {\cal F}, P )$ be a
filtered probability space. We recall that for any progressively
measurable process $\eta$ such that $E[\int_0^T|\eta_s|ds]<\infty$,
we have
$$\lim_{n\rightarrow\infty}\hat{E}[\int_0^T\delta_n(s)\eta_sds]=0.$$
Therefore, Theorem 3.3 presents a particular property of
$G$-expectation space relative to probability space.

\noindent {\bf Corollary 3.6} Let $\zeta,\eta\in M^1_G(0,T)$. If
$\int_0^t\eta_sd\langle B\rangle_s=\int_0^t\zeta_sds$ for all
$t\in[0,T]$, then
$E[\int_0^T|\eta_s|ds]=\hat{E}[\int_0^T|\zeta_s|ds]=0$.

{\bf Proof.} By Theorem 3.4, we have
$$\limsup_{n\rightarrow\infty}\hat{E}[\int_0^T\delta_n(s)\eta_sd\langle
B\rangle_s]=\lim_{n\rightarrow\infty}\hat{E}[\int_0^T\delta_n(s)\zeta_sd
s]=0.$$ By Theorem 3.3, we have $\hat{E}[\int_0^T|\eta_s|ds]=0$,
which leads to $\hat{E}[\int_0^T|\zeta_s|ds]=0$. $\Box$

The following corollary is about the uniqueness of representation
for $G$-martingales with finite variation.

\noindent {\bf Corollary 3.7} Let $\zeta,\eta\in M^1_G(0,T)$. If for
all $t\in[0,T]$,
\begin {eqnarray}\label {eqn4}\int_0^t\eta_sd\langle B\rangle_s-\int_0^t2G(\eta_s)ds=\int_0^t\zeta_sd\langle B\rangle_s-\int_0^t2G(\zeta_s)ds,
\end {eqnarray} we have $\hat{E}[\int_0^T|\eta_s-\zeta_s|ds]=0$.

{\bf Proof.} By the assumption, we have
$$\int_0^t(\eta_s-\zeta_s)d\langle
B\rangle_s=\int_0^t2[G(\eta_s)-G(\zeta_s)]ds, \ \textmd{for all} \
t\in[0,T].$$ Since $\eta-\zeta, 2[G(\eta)-G(\zeta)]\in M^1_G(0,T)$,
we have $\hat{E}[\int_0^T|\eta_s-\zeta_s|ds]=0$ by Corollary 3.6.
$\Box$

\noindent {\bf Remark 3.8}(i) In the setting considered in this
article,
$G(a)=\frac{1}{2}(\overline{\sigma}^2a^+-\underline{\sigma}^2a^-)$.
For $\varepsilon\in(0,
\frac{\overline{\sigma}^2-\underline{\sigma}^2}{2})$, [HuP10]
defined $G_\varepsilon$ in the following way:  $$G_\varepsilon(a) =
G(a) - \frac{\varepsilon}{2} |a|, \ \textmd{for all}\ a\in R.$$
Indeed, Proof to Theorem 3.3 in the next section  leads to the
following conclusion:
\begin {eqnarray}\label {eqn7}d(\eta)\geq\varepsilon
\hat{E}_{G_\varepsilon}[\int_0^T|\eta_s|ds].
\end {eqnarray}

(ii) For $\eta\in M^1_G(0,T)$, let $K_t=\int_0^T\eta_sd\langle
B\rangle_s-\int_0^T2G(\eta_s)ds$. Then, by Theorem 3.4, we have
\begin
{eqnarray}\label{eqn6}\hat{E}(-K_T)\geq\limsup_{n\rightarrow\infty}\hat{E}(\int_0^T\delta_n(s)dK_s)=d(\eta).
\end {eqnarray}
This, combined with (\ref{eqn7}), leads to the following estimate:
$$\hat{E}[-K_T]\geq\varepsilon
\hat{E}_{G_\varepsilon}[\int_0^T|\eta_s|ds],$$ which was already
proved in [HuP10]. Then for $\eta,\zeta\in M^1_G(0,T)$ such that
(\ref {eqn4}) and
\begin {eqnarray}\label
{eqn5}\int_0^t2[G(\eta_s)-G(\zeta_s)]ds=\int_0^t2[G(\eta_s-\zeta_s)]ds
\ \textmd{for all} \ t\in[0,T]
\end {eqnarray} hold, we have
$\hat{E}[\int_0^T|\eta_s-\zeta_s|ds]=0$. However, (\ref{eqn5}) does
not hold generally since the nonlinearity of $G$, which is the main
difficulty to deal with such questions. $\Box$

\section{Proof to Theorem 3.3}

In order to prove Theorem 3.3, we first introduce two lemmas.

Let $\Omega_T=C_b([0,T];R)$ be endowed with the supremum norm and
let $\sigma: [0,T]\times \Omega_T\rightarrow R$ be a measurable
mapping satisfying

i) $\sigma$ is bounded;

ii) There exists $C>0$ such that $|\sigma(s,x)-\sigma(s,y)|\leq
C\|x-y\|$ for any $s\in[0,T]$ and $x,y\in C_b([0,T];R)$;

iii)For $t\in[0,T]$, $\sigma(t,\cdot)$ is ${\cal B}_t(\Omega_T)$
measurable.

Then the following lemma is easy.

\noindent {\bf Lemma 4.1}  Let $(\Omega, F, {\cal F}, P )$ be a
filtered probability space and let $M$ be a continuous
$F$-martingale with $\langle M\rangle_t-\langle M\rangle_s\leq
C(t-s)$ for some $C>0$ and any $0\leq s<t\leq T$. Let $F^X$ be the
augmented filtration generated by $X$. Then for any $Y_0\in {\cal
F}^X_0$, there exists a unique $F$-adapted continuous process with
$E[\sup_{t\in[0,T]}|Y_t|^2]<\infty$ such that
$Y_t=Y_0+\int_0^t\sigma(s,Y)dX_s$. Moreover, $Y$ is $F^X$-adapted.
$\Box$

Let $(\Omega, {\cal F}, P )$ be a probability space and let
$\{W_t\}$ be a standard 1-dimensional Brownian motion on $(\Omega,
{\cal F}, P)$. Let $F^W$ be the augmented filtration generated by
$W$.

Denote by ${\cal A}^0([c,C])$, for some $0<c\leq C<\infty$, the
collection of $F^W$ adapted processes in the following form
$$h_s=\sum_{i=0}^{m-1}\xi_i1_{]\frac{iT}{m},\frac{(i+1)T}{m}]}(s),$$
where
$\xi_i=\psi_i(\int_{\frac{(i-1)T}{m}}^{\frac{iT}{m}}h_sdW_s,\cdot\cdot\cdot,\int_0^{\frac{T}{m}}h_sdW_s)$,
$\psi_i\in C_{b,lip}(R^i)$, $c\leq|\psi_i|\leq C$. Denote by ${\cal
A}([c,C])$ the collection of $F^W$ adapted processes such that
$c\leq|h_s|\leq C$.

\noindent {\bf Lemma 4.2} ${\cal A}^0([c,C])$ is dense in ${\cal
A}([c,C])$ under the norm $$\|h\|_2=[E(\int_0^T|h_s|^2ds)]^{1/2}.$$

{\bf Proof.} Let $h_s=\sum_{i=0}^{m-1}\xi_i1_{]{\frac{iT}{m}},
\frac{(i+1)T}{m}]}(s)$, where
$$\xi_i=\varphi_i(W_{\frac{iT}{m}}-W_{\frac{(i-1)T}{m}},\cdot\cdot\cdot,
W_{\frac{T}{m}}-W_0),$$ $$\varphi_i\in C_{b, lip}(R^i),
c\leq|\varphi_i|\leq C.$$ Then $\sigma(s,x)=h_s^{-1}(x)$ is a
bounded Lipschitz function. Let $X_t:=\int_0^th_sdW_s$. Since
$W_t=\int_0^t\sigma(s,W)dX_s$, we conclude, by Lemma 4.1, that $W$
is $F^X$-adapted.

For a process $\{X_t\}$, we denote the vector
$(X_T-X_{\frac{(m-1)T}{m}},\cdot\cdot\cdot, X_{\frac{T}{m}}-X_0)$ by
$X^m_{[0,T]}$.

For arbitrary $\varepsilon_i>0$, $i=0,\cdot\cdot\cdot,m-1$, there
exists $\psi_{i}\in C_{b,lip}(R^{in_i})$ with the Lipschitz constant
$L_{i}$ such that
$E[|\xi_{i}-\widetilde{\xi}_{i}|^2]<\varepsilon_i^2$. Here
$\widetilde{\xi}_i=\psi_i(X^{in_i}_{[0,\frac{iT}{m}]}),$
$c\leq|\psi_i|\leq C$. Without loss of generality, we assume that
there exists $K_{ji}\in N$ such that $n_j=K_{ji}n_i$ for $m-1\geq
i>j\geq0$.

Define $\widehat{\xi}_i$ in the following way:

   $\widehat{\xi}_0=\widetilde{\xi}_0$,

   For $s\in]0,\frac{T}{m}]$, $\widehat{h}_s=\widehat{\xi}_0$,

   Assume that we have defined $\widehat{h}_s$ for all $s\in[0,
   \frac{iT}{m}]$, $0\leq i\leq m-1$,

   Define $\widehat{X}_t:=\int_0^t\widehat{h}_sdW_s$, for $t\in[0,
   \frac{iT}{m}]$,

   $\widehat{\xi}_{i}=\psi_{i}(\widehat{X}^{in_i}_{[0,\frac{iT}{m}]}),$

   For $s\in]\frac{iT}{m},\frac{(i+1)T}{m}]$,
   $\widehat{h}_s=\widehat{\xi}_i$.

We claim that for any $m-1\geq i\geq1$,
\begin {eqnarray}\label{eqn3}\hat{E}[|\widehat{\xi}_i-\widetilde{\xi}_i|^2]\leq\sum_{j=0}^{i-1}A^i_j\varepsilon_j^2,
\end {eqnarray}
where $A_j^i=2TL_i^2(\sum_{k=j+1}^{i-1}A^k_j+1)$, for $i\geq j+2$,
$A^i_{i-1}=2TL_i^2$, which shows that $A_j^i$ depends only on
$L_{j+1},\cdot\cdot\cdot,L_i$ and $T$.

Indeed, $E[|\widehat{\xi}_1- \widetilde{\xi}_1|^2]\leq
L_1^2E[|\widehat{\xi}_0-
\xi_0|^2|]E[|W^{n_1}_{[0,\frac{T}{m}]}|^2]=\frac{T}{m}L_1^2\varepsilon_0^2\leq
A^1_0\varepsilon_0^2$. Assume (\ref {eqn3}) holds for $1\leq i\leq
l$. For $i=l+1$,
\begin {eqnarray*}\
& &E[|\widehat{\xi}_{l+1}-\widetilde{\xi}_{l+1}|^2]\\
&\leq&
L^2_{l+1}\sum_{i=0}^lE[|\widehat{\xi}_i-\xi_i|^2]E[|W^{n_{l+1}}_{[\frac{iT}{m},
\frac{(i+1)T}{m}]}|^2]\\
&\leq&2TL^2_{l+1}\sum_{i=0}^lE[(|\widehat{\xi}_i-\widetilde{\xi}_i|^2+|\widetilde{\xi}_i-\xi_i|^2)]\\
&\leq&2TL^2_{l+1}(\sum_{i=0}^l\varepsilon_i^2+\sum_{i=1}^l\sum_{j=0}^{i-1}A^i_j\varepsilon_j^2)\\
&=&2TL^2_{l+1}[\sum_{j=0}^{l-1}(\sum_{i=j+1}^{l}A^i_j+1)\varepsilon_j^2+\varepsilon_l^2]\\
&=&\sum_{j=0}^lA^{l+1}_j\varepsilon_j^2.
\end {eqnarray*} Then
\begin {eqnarray*}
& &E[|\widehat{\xi}_i-\xi_i|^2]\\
&\leq&2(E[|\widehat{\xi}_i-\widetilde{\xi}_i|^2]+E[|\widetilde{\xi}_i-\xi_i|^2])\\
&\leq& 2\varepsilon_i^2+2\sum_{j=0}^{i-1}A^i_j\varepsilon_j^2\\
&=:&\sum_{j=0}^{i}B^i_j\varepsilon_j^2,
\end {eqnarray*} which shows that
$B_j^i$ depends only on $L_{j+1},\cdot\cdot\cdot,L_i$ and $T$. So
for any $\varepsilon>0$, we can choose $\widehat{\xi}_i$,
$i=0,\cdot\cdot\cdot,m-1$ defined above such that
$E[|\widehat{\xi}_i-\xi_i|^2]<\varepsilon$ for all
$i=0,\cdot\cdot\cdot,m-1$. Then
$$E[\int_0^T|h_s-\widehat{h}_s|^2]<T\varepsilon.$$
$\Box$

{\bf Proof to Theorem 3.3.} For $\eta\in M^1_G(0,T)$ with
$\hat{E}[\int_0^T|\eta_s|ds]>0,$ by Theorem 2.6 and Remark 2.7,
there exists $\varepsilon>0$ and $P\in {\cal P}_M$ such that
$E_P[\int_0^T|\eta_s|ds]=:A>0$ and for any $0\leq s<t\leq T$
$$(\underline{\sigma}^2+\varepsilon)(t-s)\leq\langle B\rangle_t-\langle B\rangle_s\leq(\overline{\sigma}^2-\varepsilon)(t-s), \ \textit{P-a.s..}$$
For any
$\frac{A\varepsilon}{(\overline{\sigma}^2+\varepsilon)}>\delta>0$,
there exists $\zeta\in M^0_G(0,T)$ such that
$$\hat{E}[\int_0^T|\eta_s-\zeta_s|ds]<\delta.$$

Let $(\Omega^0, F=\{{\cal F}^0_t\}, {\cal F}^0, P^0 )$ be a filtered
probability space, and $\{W_t\}$ be a d-dimensional Brownian motion
under $P^0$.  By Remark 2.7, there exists an $F$ adapted process $h$
with $\underline{\sigma}^2+\varepsilon\leq h_s^2\leq
\overline{\sigma}^2-\varepsilon$ such that $P=P^0\circ(\int_0^\cdot
h_sdW_s)^{-1}$.

Without loss of generality, by Lemma 4.2, we assume that there
exists $m\in N$ such that
$$\zeta_s=\sum_{i=0}^{m-1}\xi_{\frac{iT}{m}}1_{]\frac{iT}{m},\frac{(i+1)T}{m}]}(s)$$
where
$\xi_{\frac{iT}{m}}=\varphi_i(B_{\frac{iT}{m}}-B_{\frac{(i-1)T}{m}},\cdot\cdot\cdot,
B_{\frac{T}{m}}-B_0)$, $\varphi_i\in C_{b, lip}(R^i)$, for all
$0\leq i\leq m-1$;
$$h_s=\sum_{i=0}^{m-1}a_{\frac{iT}{m}}1_{]\frac{iT}{m},\frac{(i+1)T}{m}]}(s)$$
where
$a_{\frac{iT}{m}}=\psi_i(\int_{\frac{(i-1)T}{m}}^{\frac{iT}{m}}h_s
dW_s,\cdot\cdot\cdot, \int_0^{\frac{T}{m}}h_s dW_s)$,
$\underline{\sigma}^2+\varepsilon\leq|\psi_i|^2\leq\overline{\sigma}^2-\varepsilon$,
$\psi_i\in C_{b, lip}(R^i)$, for all $0\leq i\leq m-1$.

1.  Define $H^i:[\underline{\sigma}^2+\varepsilon,
\overline{\sigma}^2-\varepsilon]\rightarrow[\underline{\sigma},\overline{\sigma}]$,
i=1, -1 in the following way:
$$H^1(x)^2=\overline{\sigma}^21_{[x\geq\frac{\overline{\sigma}^2+\underline{\sigma}^2}{2}]}+
(2x-\underline{\sigma}^2)1_{[x<\frac{\overline{\sigma}^2+\underline{\sigma}^2}{2}]};$$
$$H^{-1}(x)^2=(2x-\overline{\sigma}^2)1_{[x\geq\frac{\overline{\sigma}^2+\underline{\sigma}^2}{2}]}+
\underline{\sigma}^21_{[x<\frac{\overline{\sigma}^2+\underline{\sigma}^2}{2}]}.$$
It's easily seen that $H^1(x)^2+H^{-1}(x)^2=2x$ and
$H^1(x)^2-H^{-1}(x)^2\geq2\varepsilon$.

For $n\in N$, define
$H^i_n:[0,1/m]\times[\underline{\sigma}^2+\varepsilon,
\overline{\sigma}^2-\varepsilon]\rightarrow[\underline{\sigma},\overline{\sigma}]$,
$i=1,-1$ by $$H^i_n(s,x)=\sum_{j=0}^{2n-1}1_{]\frac{jT}{2mn},
   \frac{(j+1)T}{2mn}]}(s)H^{(-1)^ji}(x).$$

2. Fix $n\in N$.

   $a^n_0=a_0,$  $\xi^n_0=\xi_0,$

   For $s\in]0,\frac{T}{m}]$, $h^n_s=H^{\textmd{sgn}(\xi^n_0)}_n(s,(a^n_0)^2)$;

   Assume that we have defined $h^n_s$ for all $s\in[0,
   \frac{iT}{m}]$, $0\leq i\leq m-1$.

   $a^n_{\frac{iT}{m}}=\psi_i(\int_{\frac{(i-1)T}{m}}^{\frac{iT}{m}}h^n_s
dW_s,\cdot\cdot\cdot, \int_0^{\frac{T}{m}}h^n_s dW_s)$,

$\xi^n_{\frac{iT}{m}}=\varphi_i(\int_{\frac{(i-1)T}{m}}^{\frac{iT}{m}}h^n_s
dW_s,\cdot\cdot\cdot, \int_0^{\frac{T}{m}}h^n_s dW_s)$,

For $s\in]\frac{iT}{m},\frac{(i+1)T}{m}]$,
$h^n_s=H^{\textmd{sgn}(\xi^n_\frac{iT}{m})}_n(s-\frac{iT}{m},(a^n_{\frac{iT}{m}})^2)$.

3. $E_{P}[\int_0^T|\zeta_s|ds]=E_{P_{h^n}}[\int_0^T|\zeta_s|ds]$.

In fact,
\begin {eqnarray*}
& &E_{P}[\int_0^T|\zeta_s|ds]\\
&=&\frac{T}{m}E_{P_0}[\sum_{i=0}^{m-1}|\varphi_i(\int_{\frac{(i-1)T}{m}}^{\frac{iT}{m}}h_sdW_s,
    \cdot\cdot\cdot,\int_0^{\frac{T}{m}}h_sdW_s)|]\\
&=&:E_{P_0}[\Phi(\int_{\frac{(m-1)T}{m}}^Th_sdW_s,
    \cdot\cdot\cdot,\int_0^{\frac{T}{m}}h_sdW_s)]
\end {eqnarray*}
and
\begin {eqnarray*}
E_{P_{h^n}}[\int_0^T|\zeta_s|ds]=E_{P_0}[\Phi(\int_{\frac{(m-1)T}{m}}^Th^n_sdW_s,
    \cdot\cdot\cdot,\int_0^{\frac{T}{m}}h^n_sdW_s)].
\end {eqnarray*}

Let $x=(x_{m-1},\cdot\cdot\cdot, x_1)$. Noting that
\begin {eqnarray*}
& &\Phi_{m-1}(x)\\
&:=&E_{P_0}\{\Phi(\int_{\frac{(m-1)T}{m}}^TH^{\textmd{sgn}(\varphi_{m-1}(x))}_n(s-\frac{(m-1)T}{m},\psi_{m-1}(x)^2)dW_s,
x)\}\\
&=& E_{P_0}\{\Phi(\int_{\frac{(m-1)T}{m}}^T\psi_{m-1}(x)dW_s,x)\},
\end {eqnarray*} we have $$E_{P_0}[\Phi(\int_{\frac{(m-1)T}{m}}^Th_sdW_s,
    \cdot\cdot\cdot,\int_0^{\frac{T}{m}}h_sdW_s)]=E_{P_0}[\Phi_{m-1}(\int_{\frac{(m-2)T}{m}}^{\frac{(m-1)T}{m}}h_sdW_s,
    \cdot\cdot\cdot,\int_0^{\frac{T}{m}}h_sdW_s)]$$ $$E_{P_0}[\Phi(\int_{\frac{(m-1)T}{m}}^Th^n_sdW_s,
    \cdot\cdot\cdot,\int_0^{\frac{T}{m}}h^n_sdW_s)]=E_{P_0}[\Phi_{m-1}(\int_{\frac{(m-2)T}{m}}^{\frac{(m-1)T}{m}}h^n_sdW_s,
    \cdot\cdot\cdot,\int_0^{\frac{T}{m}}h^n_sdW_s)].$$

By induction on $m$, we get the desired result.

4. \begin {eqnarray*}
& &\hat{E}[\int_0^T\delta_{2mn}(s)\eta_sd\langle B\rangle_s]\\
&\geq& \hat{E}[\int_0^T\delta_{2mn}(s)\zeta_sd\langle
B\rangle_s]-\hat{E}[\int_0^T|\eta_s-\zeta_s|d\langle B\rangle_s]\\
&\geq& E_{P_{h^n}}[\int_0^T\delta_{2mn}(s)\zeta_sd\langle
B\rangle_s]-\overline{\sigma}^2\delta\\
&=&
E_{P_{h^n}}[\sum_{i=0}^{m-1}\int_{\frac{iT}{m}}^{\frac{(i+1)T}{m}}\delta_{2mn}(s)\zeta_sd\langle
B\rangle_s]-\overline{\sigma}^2\delta\\
&=&E_{P_{h^n}}[\sum_{i=0}^{m-1}\xi_{\frac{iT}{m}}\int_{\frac{iT}{m}}^{\frac{(i+1)T}{m}}\delta_{2mn}(s)d\langle
B\rangle_s]-\overline{\sigma}^2\delta\\
&\geq&\frac{T}{m}\varepsilon E_{P_{h^n}}[\sum_{i=0}^{m-1}|\xi_{\frac{iT}{m}}|]-\overline{\sigma}^2\delta\\
&=&\varepsilon E_{P_{h^n}}[\int_0^T|\zeta_s|ds]-\overline{\sigma}^2\delta\\
&=&\varepsilon
E_{P}[\int_0^T|\zeta_s|ds]-\overline{\sigma}^2\delta\\
&\geq&\varepsilon
E_{P}[\int_0^T|\eta_s|ds]-\varepsilon\delta-\overline{\sigma}^2\delta\\
&\geq& A\varepsilon-\varepsilon\delta-\overline{\sigma}^2\delta>0.
\end {eqnarray*} Since $A, \varepsilon, \delta$ do not depend on $n$,
we have $d(\eta)\geq
A\varepsilon-\varepsilon\delta-\overline{\sigma}^2\delta>0.$ The
proof is completed. $\Box$




\providecommand{\bysame}{\leavevmode\hbox
to3em{\hrulefill}\thinspace}
\providecommand{\MR}{\relax\ifhmode\unskip\space\fi MR }
\providecommand{\MRhref}[2]{%
  \href{http://www.ams.org/mathscinet-getitem?mr=#1}{#2}
} \providecommand{\href}[2]{#2}


\begin{thebibliography}{MOR91}
\bibitem[DHP08]{A}
Denis, L., Hu, M. and Peng S. \emph{Function spaces and capacity
related to a sublinear expectation: application to $G$-Brownian
motion pathes}. To appear in Potential Anal.
\bibitem[HP09]{A}
Hu, M. and Peng, S. (2009) \emph{On representation theorem of
Gexpectations and paths of $G$-Brownian motion}. Acta Math Appl
Sinica English Series, 25(3): 1-8.
\bibitem[HuP10]{A}
Hu, Y. and Peng, S. (2010) \emph{Some Estimates for Martingale
Representation under G-Expectation}. arXiv:1004.1098v1 [math.PR] 7
Apr 2010.

\bibitem[P07a]{A}
Peng, S. (2007) \emph{$G$-expectation, $G$-Brownian Motion and
Related Stochastic Calculus of It$\hat{o}$ type}. Stochastic
analysis and applications, 541¨C567, Abel Symp., 2, Springer,
Berlin.
\bibitem[P07b]{A}
Peng, S. (2007) \emph{$G$-Brownian Motion and Dynamic Risk Measure
under Volatility Uncertainty}. arXiv:0711.2834v1 [math.PR] 19 Nov
2007.

\bibitem[P08]{A}
Peng, S. (2008) \emph{Multi-Dimensional $G$-Brownian Motion and
Related Stochastic Calculus under $G$-Expectation}, in Stochastic
Processes and their Applications, 118(12), 2223-2253.

\bibitem[P10]{A}
Peng, S. (2010) \emph{Nonlinear Expectations and Stochastic Calculus
under Uncertainty}, arXiv:1002.4546v1 [math.PR] 24 Feb 2010.

\bibitem[Song10]{A}
Song, Y. \emph{Characterizations of processes with stationary and
independent increments under $G$-expectation} arXiv:1009.0109v1
[math.PR] 1 Sep 2010.





 \end{thebibliography}
\end{document}